\documentclass[a4paper, 10pt, conference]{ieeeconf}      

\IEEEoverridecommandlockouts                              

\usepackage{graphicx}
\usepackage[dvips]{epsfig} 
\usepackage{amsmath} 
\usepackage{amssymb}  
\usepackage{color}  
\usepackage{kotex}  
\usepackage{diagbox} 
\usepackage{siunitx} 

\title{\LARGE \bf
Minimum-Length Point-to-Line Path Planning for \\ Unmanned Aerial Vehicles
}

\author{Dongsin Kim$^{1}$ and Keumjin Lee$^{1}$ 
\thanks{$^{1}$D. Kim and K. Lee are with the School of Air Transport, Transportation and Logistics, Korea Aerospace University,
        Goyang 10540, South Korea,
        {\tt\small ds.kim@kau.kr}, {\tt\small keumjin.lee@kau.ac.kr}}%
}

\begin{document}

\maketitle
\thispagestyle{empty}
\pagestyle{empty}

\begin{abstract}
This paper presents a method of finding the shortest path for an unmanned aerial vehicle (UAV) flying from an initial point to a target line at a constant altitude. The length of a Dubins path is derived as a function of the final position on the line and then differentiated to obtain its extreme value. The primary contribution of the study is a simple analytical solution to determine the minimum-length Dubins path from an initial position to a target line with initial and final orientations given. The proposed method is demonstrated with numerical examples.
\end{abstract}

\section{INTRODUCTION}

Between two points with initial and final orientations given, Dubins showed that a minimum-length path consists of no more than three maximum-curvature circular-arcs and straight-line segments \cite{dubins1957curves}. There are six types Dubins paths between two points. One study suggested a logical classification framework to determine the shortest Dubins path directly without calculating all possible path types \cite{shkel2001classification}. Another method finds the shortest Dubins path converging tangentially to a given line \cite{hota2014optimal}, and while yet another method finds the path towards the direction tangential to a target circle \cite{manyam_shortest_2019}, \cite{chen2020dubins}. There have been various studies about possible applications of Dubins paths, including the problem of chasing a moving target \cite{manyam2018optimal}, \cite{gopalan2017generalized}, \cite{ZHENG2021109557}. Building upon previous studies, this paper provides a method to find the shortest Dubins path that starts from an initial point to a target line with initial and final orientations given.

\section{PROBLEM STATEMENT}

\begin{figure}[!htb]
    \centering
    \includegraphics[width=1\columnwidth]{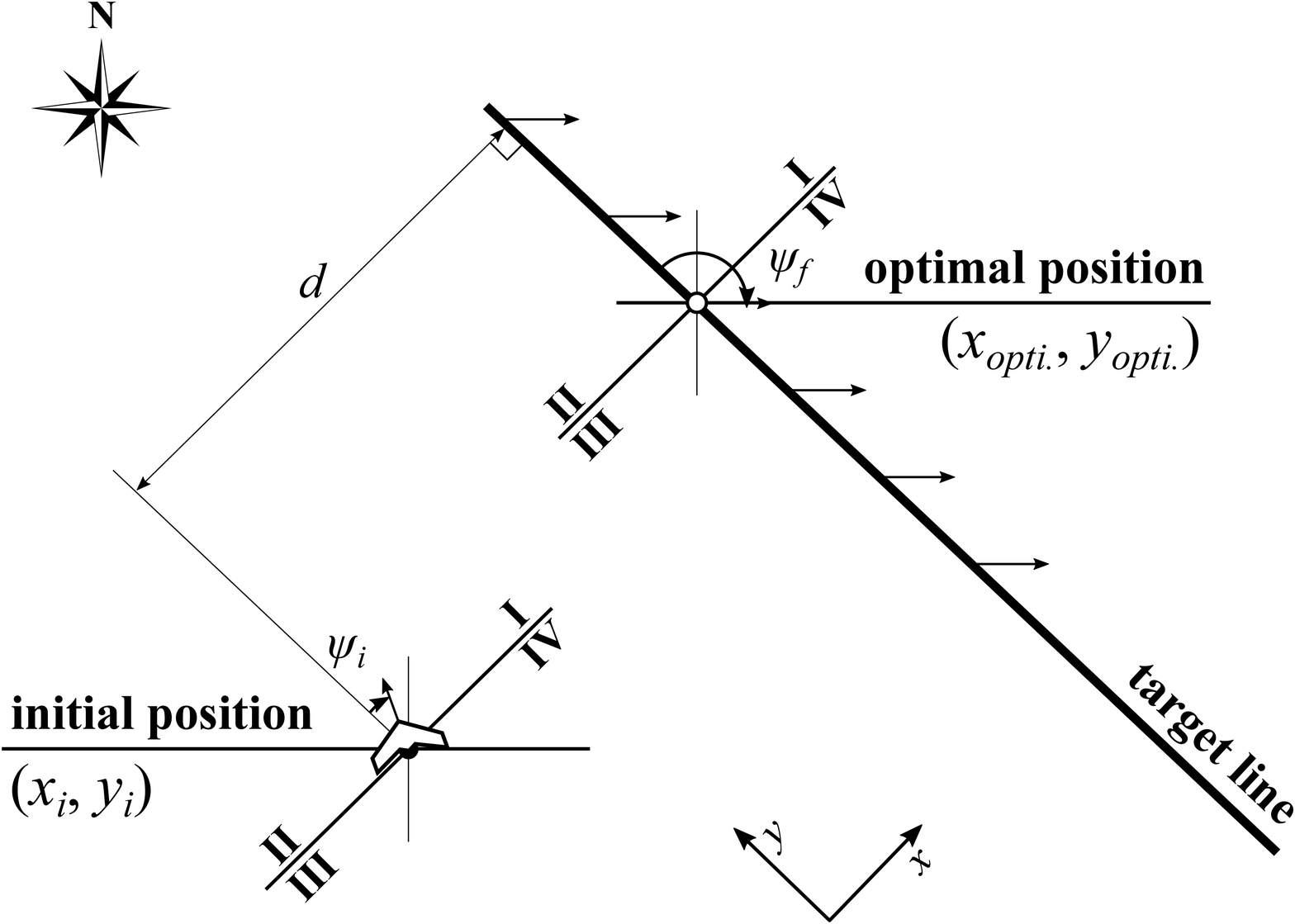}
    \caption{Notation and geometry for the point-to-line path planning}%
    \label{fig:ProblemStatement}
\end{figure}

We consider the problem of finding the shortest Dubins path from an initial position $(x_{i}, y_{i})$ to a target line in a horizontal plane, as shown in Fig.~\ref{fig:ProblemStatement}. In this study, arbitrary initial and final orientations of the vehicle are considered ($\psi_{i}$ and $\psi_{f}$, respectively). In contrast, in a previous study, it was assumed that the final orientation of the vehicle is aligned with the direction of the target line \cite{hota2014optimal}.

Between two positions with given initial and final orientations, six types Dubins paths can be generated in general \cite{dubins1957curves}. However, in this study, we assume that the perpendicular distance ($d$) between the initial position and the target line is greater than or equal to four times the vehicle’s minimum turn radius ($r$). In this case, the path type will be one of the following four combinations: RSR, RSL, LSL, or LSR, where clockwise and counter-clockwise circular arc segments are denoted by R and L, respectively, and straight-line segments are denoted by S. The results of this study can be extended to the case of $d<4r$. For brevity, we assume that the initial location of the vehicle is on the left side of the target line, but this assumption can be removed.



\section{DERIVATION OF DUBINS PATH LENGTHS}

We geometrically derived the path length from an initial point to a target line for each type of Dubins path. Without loss of generality, we assume that the initial position $(x_{i}, y_{i})$ is at the origin $(0, 0)$, and the target line is aligned with the y-axis at a distance ($d>0$) from the origin. $y_{f}$ is any point on the target line.

\subsection{RSR Path}
\begin{figure}[!htb]
    \centering
    \includegraphics[width=1\columnwidth]{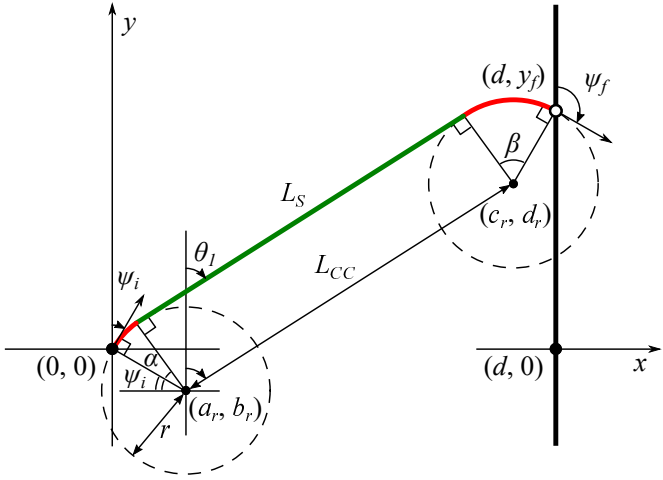}
    \caption{RSR Path}%
    \label{fig:RSR}
\end{figure}
Fig.~\ref{fig:RSR} shows an RSR path for an arbitrary position $y_{f}$. The length of the RSR path ($L_{RSR}$) is obtained by:
\begin{align}
    L_{RSR} = L_{S} + r(\alpha + \beta) \label{eqn:5}
\end{align}
where $L_{S}$ is the length of the straight-line segment, and $\alpha$ and $\beta$ represent the magnitudes of the initial and final right turns, respectively:
\begin{align}
    L_{S} = L_{CC}
\end{align}
\begin{align}
    \alpha = \theta_{1} - \psi_{i} \pmod{2\pi} \label{eqn:7}
\end{align}
\begin{align}
    \beta = \psi_{f} - \theta_{1} \pmod{2\pi} \label{eqn:8}
\end{align}
$L_{CC}$ is the distance between the centers of the two circles, and $\theta_{1}$ is the orientation of the straight-line segment relative to the y-axis.
\begin{align}
    L_{CC} = \sqrt{(c_{r}-a_{r})^2 + (d_{r}-b_{r})^2}
\end{align}
\begin{align}
    \theta_{1} = \frac{\pi}{2} - atan2\left(\frac{d_{r}-b_{r}}{c_{r}-a_{r}}\right) \pmod{\pi}
\end{align}
($a_{m}, b_{m}$) and ($c_{m}, d_{m}$) are the centers of the initial and final circular-arc segments for $m\in\{r, l\}$, which stand for right and left, respectively:
\begin{align}
    (a_{r}, b_{r}) = (r\cos\psi_{i}, - r\sin\psi_{i}) \label{eqn:1}
\end{align}    
\begin{align}
    (a_{l}, b_{l}) = (- r\cos\psi_{i}, r\sin\psi_{i}) \label{eqn:2}
\end{align}
\begin{align}
    (c_{r}, d_{r}) = (d + r\cos\psi_{f}, y_{f} - r\sin\psi_{f}) \label{eqn:3}
\end{align}
\begin{align}
    (c_{l}, d_{l}) = (d - r\cos\psi_{f}, y_{f} + r\sin\psi_{f}) \label{eqn:4}
\end{align}
The RSR path length is a piecewise continuous function with respect to $y_{f}$ due to the modulus operator, with which discontinuities occur when $\alpha$ (i.e., $\theta_{1}-\psi_{i}$) $=0$ and $\beta$ (i.e., $\psi_{f} - \theta_{1}$) $=0$.

\subsection{RSL Path}
\begin{figure}[!htb]
    \centering
    \includegraphics[width=1\columnwidth]{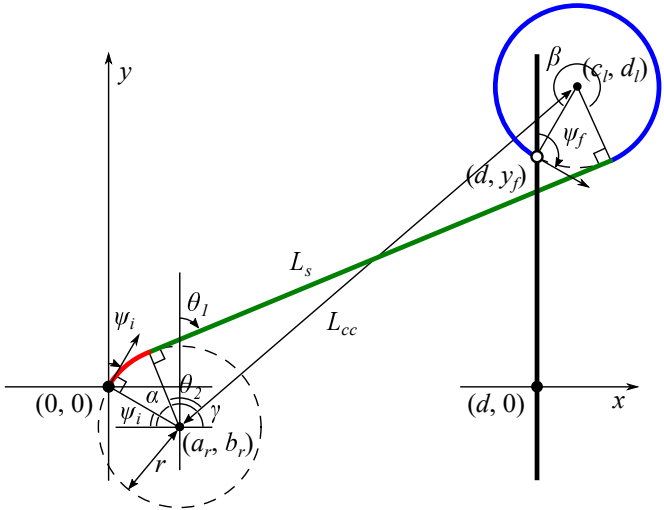}
    \caption{RSL Path}%
    \label{fig:RSL}
\end{figure}
Fig.~\ref{fig:RSL} shows an RSL path for an arbitrary position $y_{f}$. The length of the RSL path ($L_{RSL}$) is obtained by:
\begin{align}
    L_{RSL} = L_{S} + r(\alpha + \beta) \label{eqn:11}
\end{align}
where $L_{S}$ is the length of the straight-line segment, and $\alpha$ and $\beta$ represent the magnitudes of the initial right turn and the final left turn, respectively:
\begin{align}
    L_{S} = \sqrt{L_{CC}^2 - 4r^2}
\end{align}
\begin{align}
    \alpha = \theta_{1} - \psi_{i} \pmod{2\pi}
\end{align}
\begin{align}
    \beta = \theta_{1} - \psi_{f} \pmod{2\pi}
\end{align}
and
\begin{align}
    L_{CC} = \sqrt{(c_{l}-a_{r})^2 + (d_{l}-b_{r})^2}
\end{align}
\begin{align}
    \begin{split}
    \theta_{1} ={}& \frac{\pi}{2} - \left\{\gamma - \left(\frac{\pi}{2} - \theta_{2}\right)\right\}\\
    ={}& \pi - (\theta_{2} + \gamma) \pmod{\pi} \label{eqn:16}
    \end{split}
\end{align}
\begin{align}
    \theta_{2} = atan2\left(\frac{L_{S}/2}{r}\right)
\end{align}
\begin{align}
    \gamma = atan2\left(\frac{d_{l}-b_{r}}{c_{l}-a_{r}}\right)
\end{align}
$\theta_{2}$ is the interior angle of the right triangle, and $\gamma$ is the slope of the line connecting the centers of the two circles. The length of the RSL path is a piecewise continuous function with discontinuities at $\theta_{1}-\psi_{i}=0$ and $\theta_{1}-\psi_{f}=0$.

\subsection{LSL Path}
\begin{figure}[!htb]
    \centering
    \includegraphics[width=1\columnwidth]{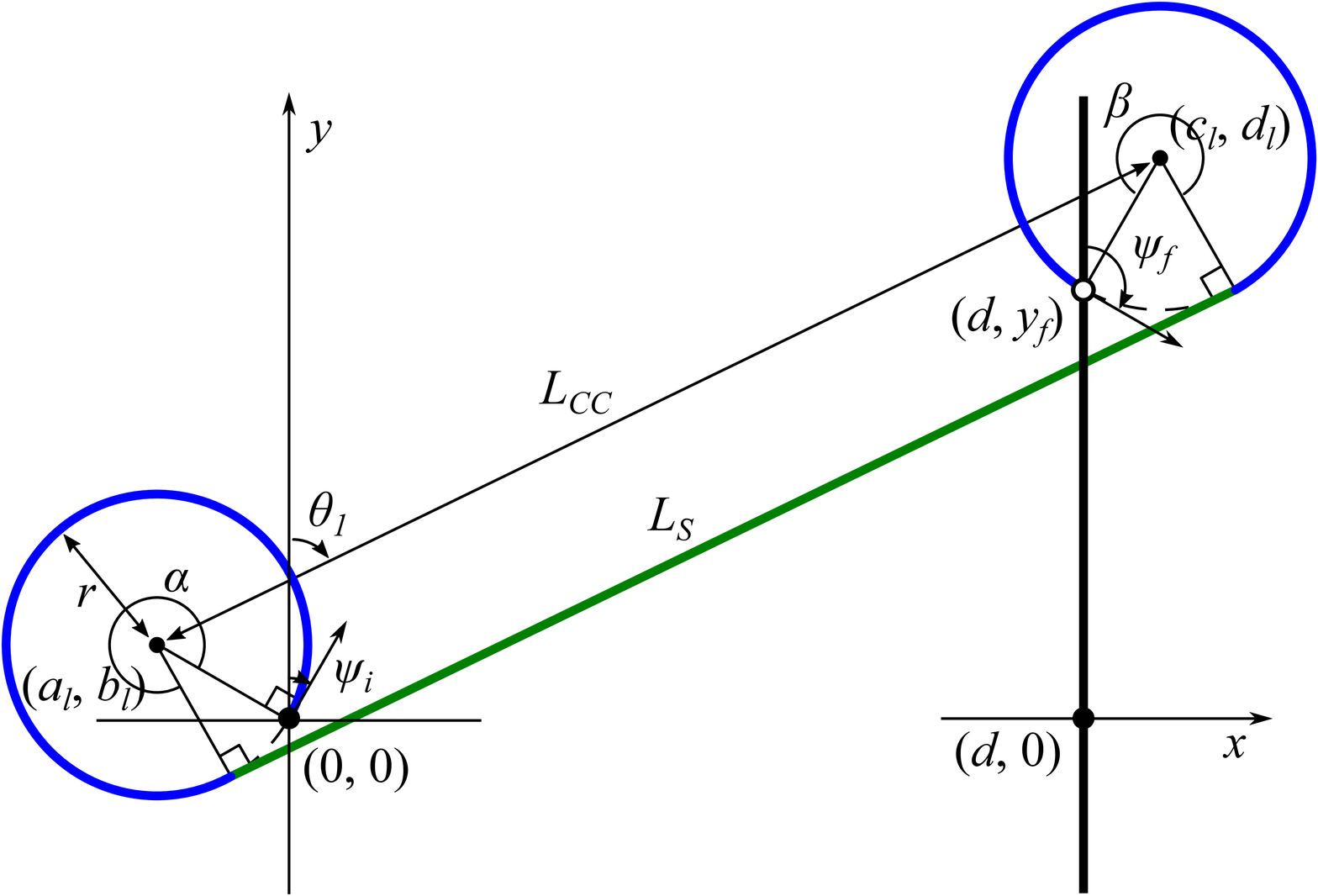}
    \caption{LSL Path}%
    \label{fig:LSL}
\end{figure}
Fig.~\ref{fig:LSL} shows an LSL path for an arbitrary position $y_{f}$. The length of the LSL path ($L_{LSL}$) is obtained by:
\begin{align}
    L_{LSL} = L_{S} + r(\alpha + \beta) \label{eqn:19}
\end{align}
where $L_{S}$ is the length of the straight-line segment, and $\alpha$ and $\beta$ represent the magnitudes of the initial and final turns, respectively:
\begin{align}
    L_{S} = L_{CC}
\end{align}
\begin{align}
    \alpha = \psi_{i} - \theta_{1} \pmod{2\pi}
\end{align}
\begin{align}
    \beta = \theta_{1} - \psi_{f} \pmod{2\pi}
\end{align}
\begin{align}
    L_{CC} = \sqrt{(c_{l}-a_{l})^2 + (d_{l}-b_{l})^2}
\end{align}
\begin{align}
    \theta_{1} = \frac{\pi}{2} - atan2\left(\frac{d_{l}-b_{l}}{c_{l}-a_{l}}\right) \pmod{\pi}
\end{align}
where $L_{LSL}$ is a piecewise continuous function with discontinuities at $\psi_{i}-\theta_{1}=0$ and $\theta_{1}-\psi_{f}=0$. 

\subsection{LSR Path}
\begin{figure}[!htb]
    \centering
    \includegraphics[width=1\columnwidth]{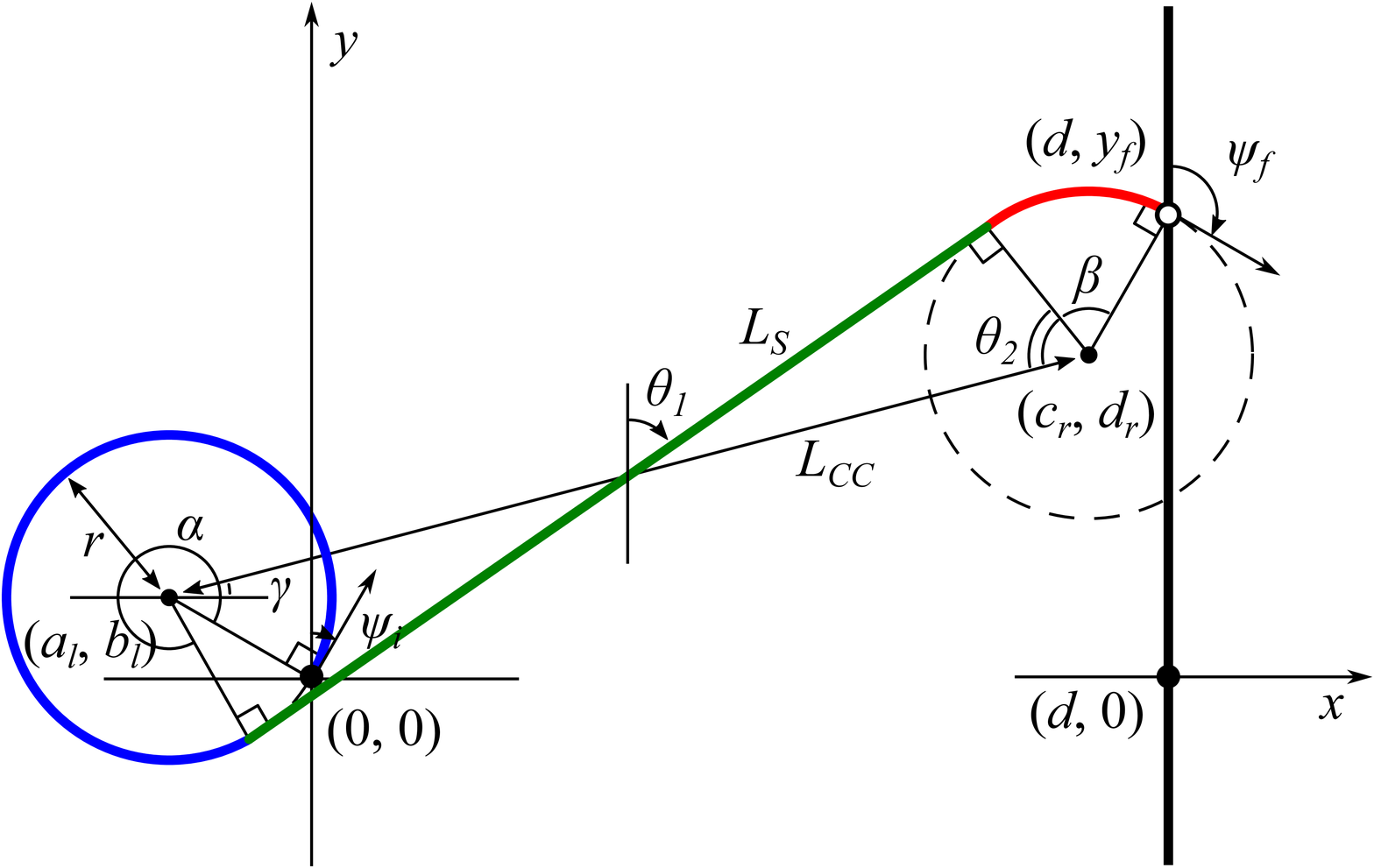}
    \caption{LSR Path}%
    \label{fig:LSR}
\end{figure}
Fig.~\ref{fig:LSR} shows an LSR path for an arbitrary position $y_{f}$. The length of the LSR path ($L_{LSR}$) is obtained by:
\begin{align}
    L_{LSR} = L_{S} + r(\alpha + \beta) \label{eqn:25}
\end{align}
where $L_{S}$ is the length of the straight-line segment, and $\alpha$ and $\beta$ represent the magnitudes of the initial and final turns, respectively:
\begin{align}
    L_{S} = \sqrt{L_{CC}^2 - 4r^2}
\end{align}
\begin{align}
    \alpha = \psi_{i} - \theta_{1} \pmod{2\pi}
\end{align}
\begin{align}
    \beta = \psi_{f} - \theta_{1} \pmod{2\pi}
\end{align}
and
\begin{align}
    L_{CC} = \sqrt{(c_{r}-a_{l})^2 + (d_{r}-b_{l})^2}
\end{align}
\begin{align}
\begin{split}
    \theta_{1} ={}& \frac{\pi}{2} - \left\{\gamma + \left(\frac{\pi}{2} - \theta_{2}\right)\right\} \label{eqn:30}\\
     ={}& \theta_{2} - \gamma \pmod{\pi}
\end{split}
\end{align}
\begin{align}
    \theta_{2} = atan2\left(\frac{L_{S}/2}{r}\right)
\end{align}
\begin{align}
    \gamma = atan2\left(\frac{d_{r}-b_{l}}{c_{r}-a_{l}}\right)
\end{align}
Discontinuities in $L_{LSR}$ occur at $\psi_{i}-\theta_{1}=0$ and $\psi_{f}-\theta_{1}=0$.

\section{OPTIMAL PATH FROM A POINT TO A LINE}

Next, we derive the optimal Dubins path from a point to a target line given initial and final orientations. To do so, the extreme path lengths for each type of Dubins path are first computed by differentiating the path length functions with respect to $y_{f}$. Then, the extreme length of one type is compared with those of the other types, and the optimal position ($y_{opti.}$) that results in the minimum length ($l_{opti.}$) from an initial point to a target line can be determined.

\subsection{Extreme Length for Each Path Type}
\subsubsection{RSR Path}
To find the extreme value of $L_{RSR}$, Eq.~\ref{eqn:5} is differentiated with respect to $y_{f}$.
\begin{align}
    \begin{split}
    \frac{d}{dy_{f}} L_{RSR} ={}& \frac{d}{dy_{f}}L_{S}\\
    ={}& \frac{d_{r}-b_{r}}{L_{S}} = 0 \label{eqn:33}
    \end{split}
\end{align}
where
\begin{align}
    \frac{d}{dy_{f}}L_{S} = \frac{d_{r}-b_{r}}{\sqrt{(c_{r}-a_{r})^2 + (d_{r}-b_{r})^2}} 
    = \frac{d_{r}-b_{r}}{L_{S}}
\end{align}

By solving Eq.~\ref{eqn:33}, we have $d_{r}-b_{r}=0$. This equation implies that the straight-line segment is normal to the target line, so we can also obtain $\theta_{1}=\pi/2$. Using Eqs.~\ref{eqn:1} and~\ref{eqn:3}, the extreme value of $L_{RSR}$ and the $y_{f}$ value for the extreme $L_{RSR}$ can be obtained as follows:
\begin{align} 
    \begin{split}
    l_{RSR}
    &= (d-r\cos\psi_{i}+r\cos\psi_{f})\\
    +& r\left(\frac{\pi}{2} - \psi_{i} \pmod{2\pi} \right)\\
    +& r\left(\psi_{f}-\frac{\pi}{2} \pmod{2\pi} \right)
    \end{split}
\end{align}
\begin{align}
    y_{f}^{RSR} = -r\sin\psi_{i} + r\sin\psi_{f}
\end{align}

\subsubsection{RSL Path}
To find the extreme value of $L_{RSL}$, Eq.~\ref{eqn:11} is differentiated with respect to $y_{f}$.
\begin{align}
    \begin{split}
    \frac{d}{dy_{f}} L_{RSL}
    ={}& \frac{d}{dy_{f}}L_{S} + 2r\frac{d}{dy_{f}}\theta_{1}\\
    ={}& \frac{d_{l}-b_{r}}{L_{S}} + 2r\left( - \frac{2r(d_{l}-b_{r})}{L_{CC}^2L_{S}} -\frac{c_{l}-a_{r}}{L_{CC}^2}\right)\\
    ={}& \frac{L_{S}(d_{l}-b_{r})-2r(c_{l}-a_{r})}{L_{CC}^2} = 0 \label{eqn:36}
    \end{split}
\end{align}
where
\begin{align}
    \frac{d}{dy_{f}}L_{S}
    = \frac{d_{l}-b_{r}}{\sqrt{(c_{l}-a_{r})^2 + (d_{l}-b_{r})^2 - 4r^2}} = \frac{d_{l}-b_{r}}{L_{S}}
\end{align}
\begin{align}
    \begin{split}
    \frac{d}{dy_{f}}\theta_{1}
    ={}& -\frac{d}{dy_{f}}\theta_{2} -\frac{d}{dy_{f}}\gamma\\
    ={}& - \frac{2r(d_{l}-b_{r})}{L_{CC}^2L_{S}} -\frac{c_{l}-a_{r}}{L_{CC}^2}
    \end{split}
\end{align}
\begin{align}
    \begin{split}
    \frac{d}{dy_{f}}\theta_{2}
    ={}& \frac{1}{1+\left(\frac{1}{2r}L_{S}\right)^2} \times \left(\frac{1}{2r}L_{S}\right)'\\
    ={}& \frac{4r^2}{L_{S}^2+4r^2} \times \left(\frac{1}{2r}\frac{d_{l}-b_{r}}{L_{S}}\right)\\
    ={}& \frac{2r(d_{l}-b_{r})}{L_{CC}^2L_{S}}
    \end{split}
\end{align}
\begin{align}
    \begin{split}
    \frac{d}{dy_{f}}\gamma
    ={}& \frac{1}{1+\left(\frac{d_{l}-b_{r}}{c_{l}-a_{r}}\right)^2} \times \left(\frac{1}{c_{l}-a_{r}}\right)\\
    ={}& \frac{c_{l}-a_{r}}{(c_{l}-a_{r})^2+(d_{l}-b_{r})^2}\\
    ={}& \frac{c_{l}-a_{r}}{L_{CC}^2}
    \end{split}
\end{align}

By solving Eq.~\ref{eqn:36}, we obtain the following equation.
\begin{align}
    \frac{d_{l}-b_{r}}{c_{l}-a_{r}} = \frac{2r}{L_{S}} \label{eqn:41}
\end{align}

If we plug Eq.~\ref{eqn:41} into Eq.~\ref{eqn:16}, we obtain:
\begin{align}
        \theta_{1} = \pi - \left\{ atan2\left(\frac{L_{S}}{2r}\right) + atan2\left(\frac{1}{\frac{L_{S}}{2r}}\right)\right\}
\end{align}

Since $L_{S}/2r>0$, we can obtain $\theta_{2}+\gamma = \pi/2$ (note that $\tan^{-1}(x)+\tan^{-1}(1/x)=\pi/2$ when $x>0$). As the straight-line segment is positioned normal to the target line, $c_{l}-a_{r}=L_{S}$. From Eq.~\ref{eqn:41}, we can also obtain $d_{l} - b_{r} = 2r$. Finally, the extreme value of $L_{RSL}$ and the $y_{f}$ value that results in the extreme $L_{RSL}$ can be obtained as follows:
\begin{align}
    \begin{split}
    l_{RSL}
    &= (d-r\cos\psi_{i}-r\cos\psi_{f})\\
    +& r\left(\frac{\pi}{2} - \psi_{i} \pmod{2\pi} \right)\\
    +& r\left(\frac{\pi}{2} - \psi_{f} \pmod{2\pi}\right)
    \end{split}
\end{align}
\begin{align}
    y_{f}^{RSL} = -r\sin\psi_{i} - r\sin\psi_{f} + 2r
\end{align}

\subsubsection{LSL Path}
To find the extreme value of $L_{LSL}$, Eq.~\ref{eqn:19} is differentiated with respect to $y_{f}$.
\begin{align}
    \frac{d}{dy_{f}} L_{LSL} = \frac{d_{l}-b_{l}}{L_{S}} = 0 \label{eqn:44}
\end{align}

By solving Eq.~\ref{eqn:44}, we have $d_{l}-b_{l}=0$, and then $\theta_{1}=\pi/2$. Using Eqs.~\ref{eqn:2} and~\ref{eqn:4}, the extreme value of $L_{LSL}$ and the $y_{f}$ value for the extreme $L_{LSL}$ can be obtained as follows:
\begin{align}
    \begin{split}
        l_{LSL} &= (d+r\cos\psi_{i}-r\cos\psi_{f})\\
        +& r\left(\psi_{i} - \frac{\pi}{2} \pmod{2\pi}\right)\\
        +& r\left(\frac{\pi}{2} - \psi_{f} \pmod{2\pi}\right)
    \end{split}
\end{align}
\begin{align}
    y_{f}^{LSL} = r\sin\psi_{i} - r\sin\psi_{f}
\end{align}

\subsubsection{LSR Path}

To find the extreme value of $L_{LSR}$, Eq.~\ref{eqn:25} is differentiated with respect to $y_{f}$.
\begin{align}
    \begin{split}
    \frac{d}{dy_{f}} L_{LSR}
    ={}& \frac{d}{dy_{f}}L_{S} - 2r\frac{d}{dy_{f}}\theta_{1}\\
    ={}& \frac{d_{r}-b_{l}}{L_{S}} - 2r\left( \frac{2r(d_{r}-b_{l})}{L_{CC}^2L_{S}} -\frac{c_{r}-a_{l}}{L_{CC}^2}\right)\\
    ={}& \frac{L_{S}(d_{r}-b_{l})+2r(c_{r}-a_{l})}{L_{CC}^2} = 0 \label{eqn:47}
    \end{split}
\end{align}
where
\begin{align}
    \begin{split}
    \frac{d}{dy_{f}}\theta_{1} ={}& \frac{d}{dy_{f}}\theta_{2} -\frac{d}{dy_{f}}\gamma\\
    ={}& \frac{2r(d_{r}-b_{l})}{L_{CC}^2L_{S}} -\frac{c_{r}-a_{l}}{L_{CC}^2}
    \end{split}
\end{align}
\begin{align}
    \begin{split}
    \frac{d}{dy_{f}}\theta_{2} 
    ={}& \frac{4r^2}{L_{S}^2+4r^2} \times \left(\frac{1}{2r}\frac{d_{r}-b_{l}}{L_{S}}\right)\\
    ={}& \frac{2r(d_{r}-b_{l})}{L_{CC}^2L_{S}}
    \end{split}
\end{align}
\begin{align}
    \begin{split}
    \frac{d}{dy_{f}}\gamma 
    ={}& \frac{c_{r}-a_{l}}{(c_{r}-a_{l})^2+(d_{r}-b_{l})^2}\\
    ={}& \frac{c_{r}-a_{l}}{L_{CC}^2}
    \end{split}
\end{align}

By solving Eq.~\ref{eqn:47}, we obtain the following equation.
\begin{align}
    \frac{d_{r}-b_{l}}{c_{r}-a_{l}} = \frac{-2r}{L_{S}} \label{eqn:52}
\end{align}

If we plug Eq.~\ref{eqn:52} into Eq.~\ref{eqn:30}, we obtain:
\begin{align}
        \theta_{1}
        = \left\{atan2\left(\frac{L_{S}}{2r}\right) + atan2\left(\frac{1}{\frac{L_{S}}{2r}}\right)\right\}
\end{align}

Since $L_{S}/2r>0$, we can obtain $\theta_{1}=\pi/2$ (note that $\tan^{-1}(-x)$ = $-\tan^{-1}(x)$), and then, $c_{r}-a_{l}=L_{S}$. From Eq.~\ref{eqn:52}, we can also obtain $d_{r} - b_{l} = -2r$. Finally, the extreme value of $L_{LSR}$ and the $y_{f}$ value that results in the extreme $L_{LSR}$ can be obtained as follows:
\begin{align}
    \begin{split}
    l_{LSR} &= (d+r\cos\psi_{i}+r\cos\psi_{f})\\
    +& r\left(\psi_{i} - \frac{\pi}{2} \pmod{2\pi}\right)\\
    +& r\left(\psi_{f} - \frac{\pi}{2}\pmod{2\pi}\right)
    \end{split}
\end{align}
\begin{align}
    y_{f}^{LSR} = r\sin\psi_{i} + r\sin\psi_{f} - 2r
\end{align}

\clearpage

\subsection{Optimal Dubins Path from a Point to a Line}

The extreme lengths for all path types and the corresponding $y_{f}$ values obtained thus far can be summarized as follows:
\begin{align} 
    \begin{split}
    l_{RSR} {}&= (d-r\cos\psi_{i}+r\cos\psi_{f})\\
    +& r\left(\frac{\pi}{2} - \psi_{i} \pmod{2\pi} \right)\\
    +& r\left(\psi_{f}-\frac{\pi}{2} \pmod{2\pi} \right)\\
    l_{RSL} {}&= (d-r\cos\psi_{i}-r\cos\psi_{f})\\
    +& r\left(\frac{\pi}{2} - \psi_{i} \pmod{2\pi} \right)\\
    +& r\left(\frac{\pi}{2} - \psi_{f} \pmod{2\pi}\right)\\
    l_{LSL} {}&= (d+r\cos\psi_{i}-r\cos\psi_{f})\\
    +& r\left(\psi_{i} - \frac{\pi}{2} \pmod{2\pi}\right)\\
    +& r\left(\frac{\pi}{2} - \psi_{f} \pmod{2\pi}\right)\\
    l_{LSR} {}&= (d+r\cos\psi_{i}+r\cos\psi_{f})\\
    +& r\left(\psi_{i} - \frac{\pi}{2} \pmod{2\pi}\right)\\
    +& r\left(\psi_{f} - \frac{\pi}{2}\pmod{2\pi}\right)
    \end{split}
\end{align}
\begin{align}
    \begin{split}
    y_{f}^{RSR} ={}& -r\sin\psi_{i} + r\sin\psi_{f}\\
    y_{f}^{RSL} ={}& -r\sin\psi_{i} - r\sin\psi_{f} + 2r\\
    y_{f}^{LSL} ={}& r\sin\psi_{i} - r\sin\psi_{f}\\
    y_{f}^{LSR} ={}& r\sin\psi_{i} + r\sin\psi_{f} - 2r
    \end{split}
\end{align}
As we seek the optimal path from an initial point to a target line, the optimal path type that results in the minimum length can be determined by comparing the length of Dubins paths as follows:
\begin{align}
    l_{opti.} = \min\{l_{RSR}, \hspace{0.2cm} l_{RSL}, \hspace{0.2cm} l_{LSL}, \hspace{0.2cm} l_{LSR}\} \label{eqn:OP}
\end{align}
\begin{align}
    y_{opti.} \in \{y_{f}^{RSR}, \hspace{0.2cm} y_{f}^{RSL}, \hspace{0.2cm} y_{f}^{LSL}, \hspace{0.2cm} y_{f}^{LSR}\}
\end{align}

\subsection{Optimal Path Decision Table}

We propose that it is not necessary to solve Eq.~\ref{eqn:OP} to find the optimal path when $\psi_{i}$ and $\psi_{f}$ are given. It is enough to consider only one type of Dubins path. The optimal path types for all possible combinations of initial and final orientations are listed in Table~\ref{tb:1}.
\begin{table}[!htb]
\caption{The Decision Table for the Optimal Path}
\label{tb:1}
\begin{center}
\begin{tabular}{c|cc|cc}
\diagbox[width=10em]{Initial \\Quadrant}{Final \\Quadrant} & I & II & III & IV\\
\hline
  &     &     &     &\\
I & RSL & RSL & RSR & RSR\\
  &     &     &     &\\
  &     &     &     &\\
II & RSL & RSL & RSR & RSR\\
  &     &     &     &\\
\hline
  &     &     &     &\\
III & LSL & LSL & LSR & LSR\\
  &     &     &     &\\
  &     &     &     &\\
IV  & LSL & LSL & LSR & LSR\\
  &     &     &     &\\
\end{tabular}
\end{center}
\end{table}

\section{EXAMPLES}

\begin{figure}[!htb]
    \centering
    \includegraphics[width=1\columnwidth]{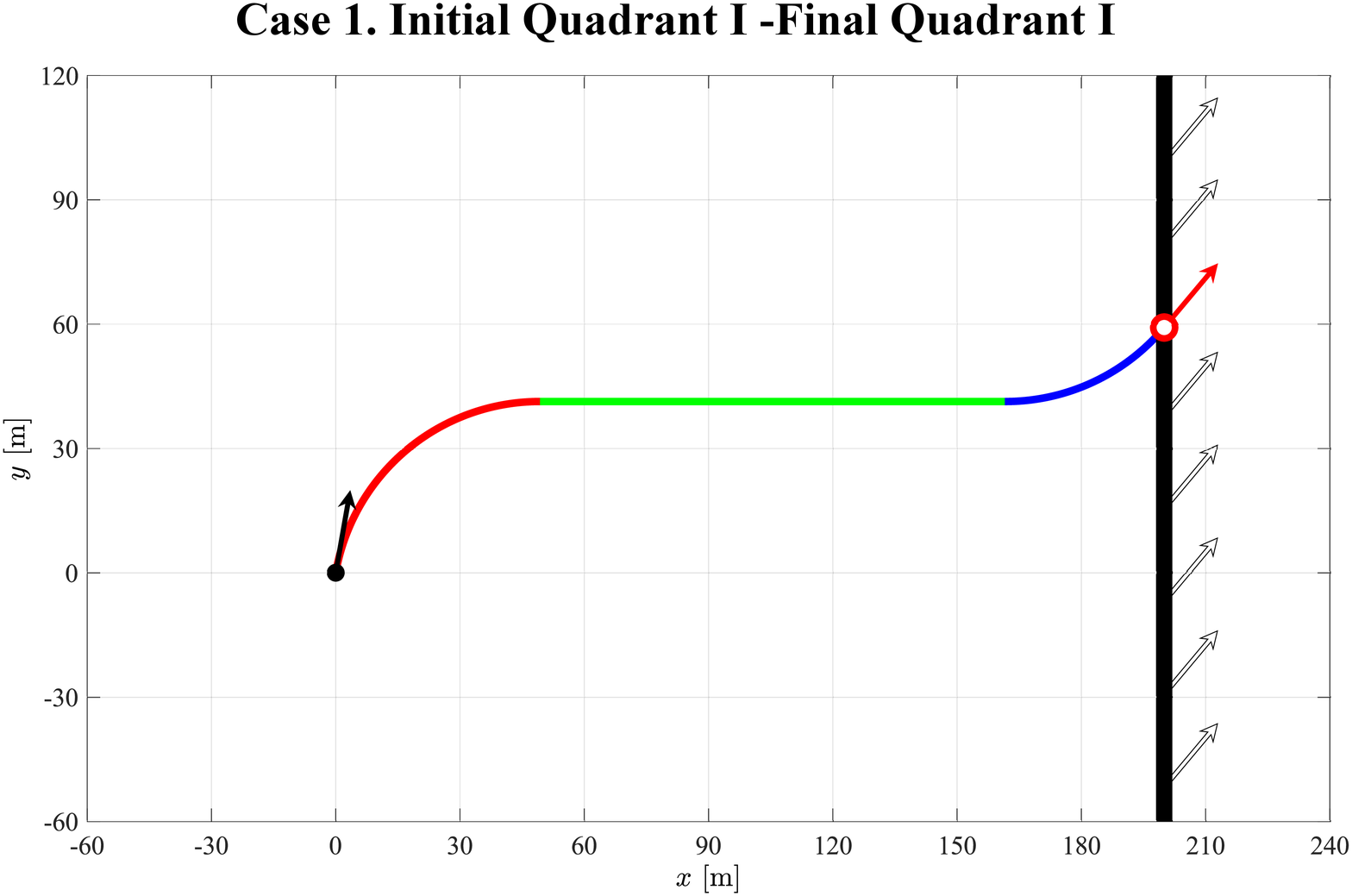}
    \caption{The optimal Dubins path found by the point-to-line path planner for initial conditions $x_{i}=0 \text{ m}$, $y_{i}=0 \text{ m}$, $\psi_{i}=\ang{10}$, $\psi_{f} = \ang{40}$, $d = 200 \text{ m}$, and $r = 50 \text{ m}$. The optimal path is RSL path with $y_{opti.} = 59.1782 \text{ m}$ and $l_{opti.} = 225.9038 \text{ m}$. Note that the slope of the straight-line segment is normal to the target line.}
    \label{fig:6}
\end{figure}

\begin{figure}[!htb]
    \centering
    \includegraphics[width=1\columnwidth]{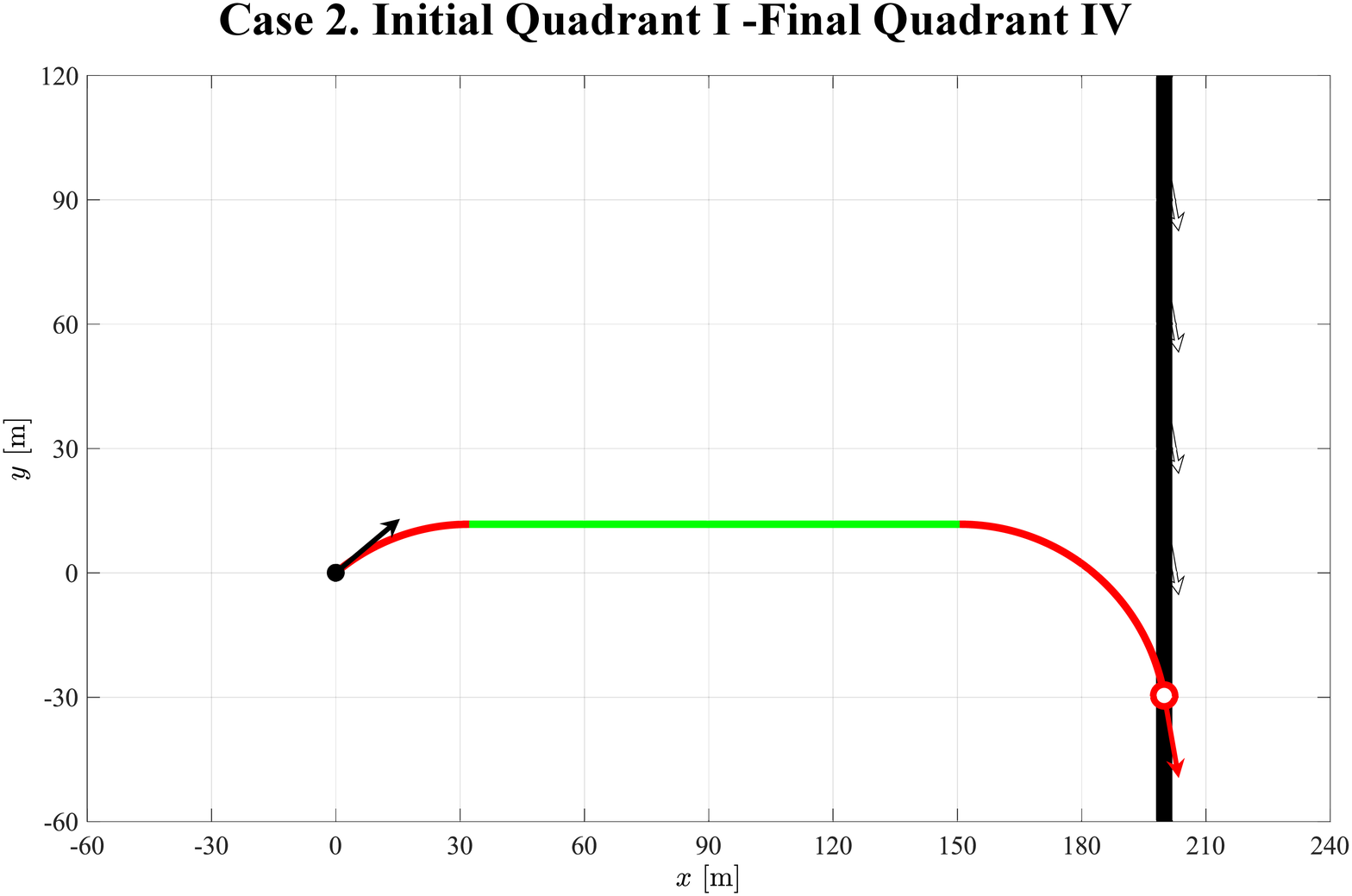}
    \caption{The optimal Dubins path found by the point-to-line path planner for initial conditions $x_{i}=0 \text{ m}$, $y_{i}=0 \text{ m}$, $\psi_{i}=\ang{50}$, $\psi_{f} = \ang{170}$, $d = 200 \text{ m}$, and $r = 50 \text{ m}$. The optimal path is RSR path with $y_{opti.} = -29.6198 \text{ m}$ and $l_{opti.} = 223.3400 \text{ m}$. Note that the slope of the straight-line segment is normal to the target line.}
    \label{fig:7}
\end{figure}

\begin{figure}[!htb]
    \centering
    \includegraphics[width=1\columnwidth]{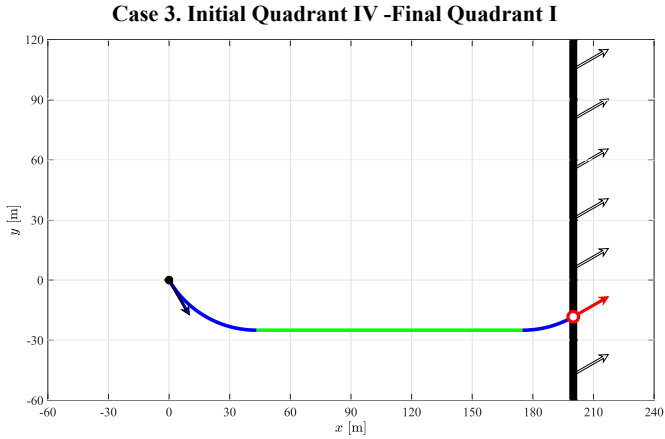}
    \caption{The optimal Dubins path found by the point-to-line path planner for initial conditions $x_{i}=0 \text{ m}$, $y_{i}=0 \text{ m}$, $\psi_{i}=\ang{150}$, $\psi_{f} = \ang{60}$, $d = 200 \text{ m}$, and $r = 50 \text{ m}$. The optimal path is LSL path with $y_{opti.} = -18.3013 \text{ m}$ and $l_{opti.} = 210.2385 \text{ m}$. Note that the slope of the straight-line segment is normal to the target line.}
    \label{fig:8}
\end{figure}

\begin{figure}[!htb]
    \centering
    \includegraphics[width=1\columnwidth]{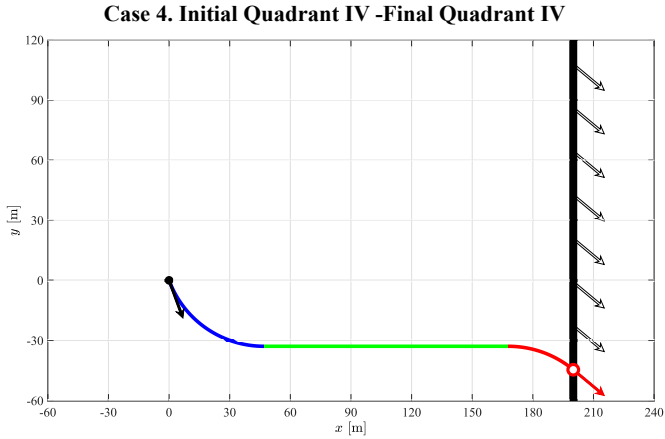}
    \caption{The optimal Dubins path found by the point-to-line path planner for initial conditions $x_{i}=0 \text{ m}$, $y_{i}=0 \text{ m}$, $\psi_{i}=\ang{160}$, $\psi_{f} = \ang{130}$, $d = 200 \text{ m}$, and $r = 50 \text{ m}$. The optimal path is LSR path with $y_{opti.} = -44.5968 \text{ m}$ and $l_{opti.} = 216.8691 \text{ m}$. Note that the slope of the straight-line segment is normal to the target line.}
    \label{fig:9}
\end{figure}

As shown in Figs.~\ref{fig:6} to~\ref{fig:9}, the method described is demonstrated with variables $x_{i}$, $y_{i}$, $\psi_{i}$, $\psi_{f}$, $d$, and $r$. In the examples, only two variables ($\psi_{i}$ and $\psi_{f}$) representing one of the four quadrants were changed in each case.

\section{CONCLUSION}

We have presented a simple geometric method to determine the shortest path from an initial position to a target line based on Dubin’s classical results for the minimum-length path between two positions with initial and final orientations. The results of the method are expressed in analytical formulae and would thus be easy to implement in real-time on-board processing for autonomous vehicles. We are currently applying the proposed method to solve the path-planning problem for a vehicle flying from a point to a region at a constant altitude, but other types of use cases should be further identified as future work.

\addtolength{\textheight}{-12cm}   



\bibliographystyle{ieeetr}
\bibliography{refs.bib}

\end{document}